\documentclass[12pt]{article}

\usepackage{amssymb,amsmath}

\begin{document}

\def\s{\subseteq}
\def\h{\widehat}
\def\v{\varphi}
\def\t{\widetilde}
\def\ov{\overline}
\def\L{\Lambda}
\def\l{\lambda}
\def\O{\Omega}
\def\H{I\!\! H}
\def\a{\approx}
\def\k{\widetilde}
\def\la{\lambda}
\def\d{\delta}
\def\L{\Lambda}
\def\O{\Omega}
\def\r{\rho}
\def\ov{\overline}
\def\un{\underline}

\title{{\bf A study concerning splitting and
jointly continuous topologies on $C(Y,Z)$}}

\author{Dimitris Georgiou$^a$, Athanasios Megaritis$^b$,\\ Kyriakos Papadopoulos$^c$, and Vasilios
Petropoulos$^a$\\
\small{$^a$University of Patras, Department of Mathematics, Greece}\\
\small{$^b$Technological Educational Institute of Western Greece,}\\[-0.8ex]
\small{Department of Accounting and Finance, 302 00 Messolonghi, Greece}\\
\small{$^c$ Section of Mathematics, College of Engineering,}\\[-0.8ex]
\small{American University of the Middle East, 15453 Dasman, Kuwait}}

\date{}
\maketitle

\begin{abstract}
Let $Y$  and $Z$ be two fixed topological spaces and $C(Y,Z)$ the
set of all continuous maps from $Y$ into $Z$. We
construct and study topologies on $C(Y,Z)$ that we call ${\cal F}_n(\tau_n)$-family-open
topologies. Furthermore, we find necessary and sufficient conditions such that these topologies  to be splitting and jointly continuous. Finally, we present questions concerning a further study on this area.
\end{abstract}

\medskip
\noindent {\bf Keywords:} Function spaces, Splitting topology,
Jointly continuous topology.

\medskip
\noindent
{\bf 2010 AMS Subject Classification.} 54C35

\newtheorem{definition}{Definition}[section]
\newtheorem{remark}{Remark}[section]
\newtheorem{remarks}{Remarks}[section]
\newtheorem{notation}{Notation}[section]
\newtheorem{examples}{Examples}[section]
\newtheorem{example}{Example}[section]
\newtheorem{theorem}{Theorem}[section]
\newtheorem{pro}{Proposition}[section]
\newtheorem{co}{Corollary}[section]
\newtheorem{lemma}{Lemma}[section]

\section{Introduction}

By $Y$ and $Z$ we denote two fixed topological spaces and
by $C(Y,Z)$ the set of all continuous maps from $Y$ to $Z$. If $t$
is a topology on $C(Y,Z)$, then the corresponding topological space
is denoted by $C_{t}(Y,Z)$.

\medskip
\noindent Let $X$ be a space, $F:X\times Y\to Z$ a continuous map
and $x\in X$. By $F_x$ we denote the continuous map from $Y$ to $Z$,
for which $F_x(y)=F(x,y)$ for every $y\in Y$. Also, by $\h{F}$ we
denote the map from $X$ to $C(Y,Z)$, for which $\h{F}(x)=F_x$, for
every $x\in X$.

\medskip
\noindent
Let $G$ be a map from $X$ to $C(Y,Z)$. By $\t{G}$ we
denote the map from $X\times Y$ to $Z$, for which
$\t{G}(x,y)=G(x)(y)$ for every $(x,y)\in X\times Y$.

\medskip
\noindent
A topology $t$ on $C(Y,Z)$ is called {\it splitting}, if
for every space $X$, the continuity of the map $F:X\times Y\to Z$
implies the continuity of the map $\h{F}:X\to C_{t}(Y,Z)$. A
topology $t$ on $C(Y,Z)$ is called  {\it jointly continuous}, if for
every space $X$, the continuity of the map $G:X\to C_{t}(Y,Z)$
implies the continuity of the map $\t{G}:X\times Y\to Z$ (see
\cite{ARENS} and \cite{AREDUG}). Let ${\cal A}$ be a fixed family of
topological spaces. If, in the above definitions, we assume that
the space $X$ belongs to ${\cal A}$, then the topology $t$  is
called {\it ${\cal A}$-splitting} (respectively,  {\it ${\cal
A}$-jointly continuous}) (see \cite{GEORGIOU1}).

\medskip
\noindent 
The Scott topology $\tau_{Sc}(P)$ on a poset $(P,\leq)$ (see, for example, \cite{SCOT}) is the family of all subsets $\H$ of $P$ such that:\\
(a) $\H=\,\uparrow\!\H$, where $\uparrow\!\H=\bigcup_{x\in\H}\{y\in P:x\leq y\}$.\\
(b) For every directed subset $D$ of $P$,
$\sup(D)\in\H$ implies that $D\cap\H\neq\emptyset$.

\medskip
\noindent 
Let $Y$ be a topological space and let ${\cal O}(Y)$ be
the family of all open subsets of $Y$ ordered via inclusion.
Then, the Scott topology $\tau_{Sc}({\cal O}(Y))$ on ${\cal O}(Y)$ is the family of all subsets $\H$ of ${\cal O}(Y)$ such that:\\
(a) $U\in\H$, $V\in{\cal O}(Y)$ and $U\subseteq V$ imply that $V\in\H$.\\
(b) For every family $\{U_i:i\in I\}\subseteq{\cal O}(Y)$ such that
$\bigcup\{U_i:i\in I\}\in\H$, there exists a finite subset $J$ of
$I$ such that $\bigcup\{U_i:i\in J\}\in\H$.

\medskip
\noindent
The {\it Isbell topology} on $C(Y,Z)$ (see, for example, \cite{MN}),
denoted here by $t_{Is}$ is the topology with subbasis: $$(\H, U)=\{f\in C(Y,Z):f^{-1}(U)\in \H\},$$ where $\H\in\tau_{Sc}({\cal O}(Y))$ and $U\in {\cal O}(Z)$.

\medskip
\noindent A subset $K$ of a topological space $Y$ is said to be {\it
relatively compact} (see, for example, \cite{LAMPAP}) if every open cover of
$Y$ has a finite subcover for $K$.

\medskip
\noindent A space $Y$ is called {\it corecompact} (see, for example,
\cite{SCOT}) if for every $y\in Y$ and for every open neighborhood
$U$ of $y$, there exists an open neighborhood $V$ of $y$, such that
$V$ is relatively compact in the space $U$.

\medskip
\noindent
Let ${\bf S}$ be the Sierpi\'{n}ski space, that is, ${\bf S}=\{0,1\}$ with the topology $\{\emptyset,\{1\},{\bf S}\}$. If $Y$ is another topological space, then $C(Y,{\bf S})=\{{\mathcal X}_{U}:U\in{\cal O}(Y)\}$, where ${\mathcal X}_U:Y\to{\bf S}$ denotes the
characteristic function of $U$, $${\mathcal X}_U(y)=\begin{cases}1 & {\rm if} \ y\in U,\\ 0 &
{\rm if} \ y\in Y\setminus U.\end{cases}$$

\medskip
\noindent Below we give some known results:

\medskip
\noindent (1) The Isbell topology on
$C(Y,Z)$ is always splitting (see \cite{MN}).

\medskip
\noindent
(2) The Isbell topology on $C(Y,Z)$ is jointly continuous
if $Y$ is a corecompact space. In this case the Isbell topology is
the greatest splitting topology (see, for example, \cite{MN} and
\cite{SCOT}).

\medskip
\noindent 
(3) Each splitting topology is contained in each jointly
continuous topology (see \cite{DUG}).

\medskip
\noindent 
(4) A topology which is larger than a jointly continuous
topology is also jointly continuous (see \cite{DUG}).

\medskip
\noindent 
(5) A topology which is smaller than a splitting topology
is also splitting (see \cite{DUG}).

\medskip
\noindent
In this paper we construct and then study topologies on the function space
$C(Y,Z)$, that we call ${\cal F}_n(\tau_n)$-family-open topologies. Furthermore, we find
necessary and sufficient conditions, for these topologies to be
splitting and jointly continuous. Section 1 contains some preliminary definitions and results. In section 2 we define topologies on $C(Y,Z)$ that we call ${\cal F}_n(\tau_n)$-family-open topologies on the set $C(Y,Z)$. In section 3 we give basic remarks for the ${\cal F}_n(\tau_n)$-family-open topologies. In sections 4 and 5 we give a characterization of ${\cal A}$-splitting and ${\cal A}$-jointly continuous topologies for  the ${\cal
F}_n(\tau_n)$-family-open topologies. Finally, in section 6 we present questions concerning the ${\cal F}_n(\tau_n)$-family-open topologies.

\section{${\cal F}_n(\tau_n)$-family-open topologies on the set\\ $C(Y,Z)$}

In what follows, the power set of a set $X$ will be denoted by ${\cal P}(X)$. 

\begin{definition}{\rm
Let $Y$ and $Z$ be two topological spaces and ${\cal O}(Y)$
the family of all open subsets of $Y$. We define topologies on the set $C(Y,Z)$ using the following steps:

\bigskip
\noindent
{\bf Step 1.} Let ${\cal F}_1\subseteq{\cal P}({\cal O}(Y))$. For every open set $U$ in the space $Y$, we set $$O^1(U)=\{\varphi\in {\cal F}_1:U\in\varphi\}.$$

\medskip
\noindent
We denote by ${\cal O}^1({\cal F}_1)$ the set $$\{O^1(U):U\in {\cal O}(Y)\}.$$ Let $\tau_1$ be a topology on ${\cal O}^1({\cal F}_1)$. By
$t_{{\cal F}_1}(\tau_1)$ we denote the topology on $C(Y,Z)$ for which the
family of all  sets of the form $$<I\!\! H, U>_1=\{f\in C(Y,Z):O^1(f^{-1}(U))\in I\!\! H\}$$ forms a subbasis, where $I\!\! H\in\tau_1$ and $U\in {\cal
O}(Z)$.

\medskip
\noindent
The topology $t_{{\cal F}_1}(\tau_1)$ on $C(Y,Z)$ will be called
{\it ${\cal F}_1(\tau_1)$-family-open topology}.

\bigskip
\noindent
{\bf Step 2.} Let ${\cal F}_2\subseteq{\cal P}({\cal O}^1({\cal
F}_1))$. For every open set $U$ in the space $Y$, we set $$O^2(U)=\{\varphi\in{\cal F}_2:O^1(U)\in\varphi\}.$$ We denote by ${\cal O}^2({\cal F}_2)$ the set $$\{O^2(U):U\in {\cal O}(Y)\}.$$

\medskip
\noindent
Let $\tau_2$ be a topology on ${\cal O}^2({\cal F}_2)$. By
$t_{{\cal F}_2}(\tau_2)$ we denote the topology on $C(Y,Z)$ for which the
family of all  sets of the form $$<I\!\! H, U>_2=\{f\in C(Y,Z):O^2(f^{-1}(U))\in I\!\! H\}$$ forms a subbasis, where $I\!\! H\in\tau_2$ and $U\in {\cal
O}(Z)$.

\medskip
\noindent
The topology $t_{{\cal F}_2}(\tau_2)$ on $C(Y,Z)$ will be called
{\it ${\cal F}_2(\tau_2)$-family-open topology}.

\bigskip
\noindent
We continue in the same manner to a Step 3, Step 4, etc.
Step $n$ will look as follows.

\bigskip
\noindent
{\bf Step {\textit n.}} Let ${\cal F}_n\subseteq{\cal P}({\cal O}^{n-1}({\cal
F}_{n-1}))$. For every open set $U$ in the space $Y$, we set $$O^n(U)=\{\varphi\in{\cal F}_{n}:O^{n-1}(U)\in\varphi\}.$$ In addition, we denote by ${\cal O}^n({\cal F}_n)$ the
set $$\{O^n(U):U\in {\cal O}(Y)\}.$$ Let $\tau_n$ be a topology on ${\cal O}^n({\cal F}_n)$. By
$t_{{\cal F}_n}(\tau_n)$ we denote the topology on $C(Y,Z)$ for which the
family of all sets of the form $$<I\!\! H, U>_n=\{f\in C(Y,Z):O^n(f^{-1}(U))\in I\!\! H\}$$ forms a subbasis, where $I\!\! H\in\tau_n$ and $U\in {\cal
O}(Z)$.

\medskip
\noindent
The topology $t_{{\cal F}_n}(\tau_n)$ on $C(Y,Z)$ will be called
{\it ${\cal F}_n(\tau_n)$-family-open topology}.}
\end{definition}

\begin{example}
{\rm {\bf (1)} Let ${\bf S}$ be the Sierpi\'{n}ski space, that is, ${\bf S}=\{0,1\}$ with the topology $\tau=\{\emptyset,\{1\},{\bf S}\}$ and let $Y=\{a,b,c\}$ with the topology $${\cal O}(Y)=\{\emptyset,U_1=\{a\},U_2=\{b\},U_{3}=\{a,b\},Y\}.$$ The family $$\tau_0=\{\emptyset,\{U_3\},\{U_1,U_3\},\{U_2,U_3\},\{U_1,U_2,U_3\},{\cal O}(Y)\}$$ defines a topology on ${\cal O}(Y)$. We consider the
topology $t(\tau_0)$ on $C(Y,{\bf S})$ which has as subbasis the family of all
sets of the form $$<I\!\! H, U>=\{f\in C(Y,{\bf S}):f^{-1}(U)\in I\!\!
H\},$$ where $I\!\! H\in\tau_0$ and $U\in\{\emptyset,\{1\},\{0,1\}\}$. 

\bigskip

$\bullet$ Let ${\cal F}_1=\{{\cal O}(Y)\}\subseteq{\cal P}({\cal O}(Y))$. We have
$$O^1(\emptyset)=O^1(U_1)=O^1(U_2)=O^1(U_3)=O^1(Y)=\{{\cal O}(Y)\}.$$
Therefore, ${\cal O}^1({\cal F}_1)=\{\{{\cal O}(Y)\}\}=\{{\cal F}_1\}$. We then consider the
topology $\tau_1=\{\emptyset,\{{\cal F}_1\}\}$ on ${\cal O}^1({\cal
F}_1)$ and define, using the {\bf Step 1},
the topology $t_{{\cal F}_{1}}(\tau_1)$ on the set $C(Y,{\bf S})$ as follows:

\medskip
\noindent
The ${\cal F}_1(\tau_1)$-family-open topology $t_{{\cal F}_1}(\tau_1)$ on
$C(Y,{\bf S})$ has as subbasis the family of all sets of the form
$$<I\!\! H, U>_1=\{f\in C(Y,{\bf S}):f^{-1}(U)\in I\!\! H\},$$ where
$I\!\! H\in\tau_1$ and $U\in\{\emptyset,\{1\},\{0,1\}\}$.

\medskip
\noindent We observe that $t_{{\cal F}_{1}}(\tau_1)\neq t(\tau_0)$ on
$C(Y,{\bf S})$. Indeed, let ${I\!\! H}_0=\{U_3\}\in\tau_0$. Then,
$<{I\!\! H}_0,\{1\}>=\{{\mathcal X}_{U_3}\}\in t(\tau_0)$. We prove that $<{I\!\! H}_0,\{1\}>\notin t_{{\cal
F}_{1}}(\tau_1)$. Indeed, we observe that
\begin{eqnarray*}
<\{{\cal F}_1\},\{1\}>_{1}&=&\{f\in C(Y,{\bf S}):O^1(f^{-1}(\{1\}))\in\{\{{\cal O}(Y)\}\}\}\\
&=&\{f\in C(Y,{\bf S}):O^1(f^{-1}(\{1\}))=\{{\cal O}(Y)\}\}\\
&=&\{{\mathcal X}_{\emptyset},{\mathcal X}_{U_1},{\mathcal X}_{U_2},{\mathcal
X}_{U_3},{\mathcal X}_{Y}\}=C(Y,{\bf S}).
\end{eqnarray*}
Thus, $t_{{\cal F}_{1}}(\tau_1)=\{\emptyset,C(Y,{\bf S})\}$ and, therefore,
$t_{{\cal F}_{1}}(\tau_1)\neq t(\tau_0)$.

\bigskip

$\bullet$ Now, let ${\cal F}'_1=\{\{U_1,U_3\},\{U_2,U_3\}\}$. Then, $${\cal O}^1({\cal F}'_1)=\{\emptyset,{\cal F}'_1,\{\{U_1,U_3\}\},\{\{U_2,U_3\}\}\}.$$ We consider the topologies $$\tau'_1=\{\emptyset,\{{\cal F}'_1,\{\{U_1,U_3\}\},{\cal O}^1({\cal F}'_1)\}\}$$ and $$\tau''_1=\{\emptyset,\{{\cal F}'_1,\{\{U_2,U_3\}\},{\cal O}^1({\cal F}'_1)\}\}$$ on ${\cal O}^1({\cal
F}'_1)$ and define, using the {\bf Step 1},
the topologies $t_{{\cal F}'_{1}}(\tau'_1)$ and $t_{{\cal F}'_{1}}(\tau''_1)$ on the set $C(Y,{\bf S})$. We observe that $$t_{{\cal F}'_{1}}(\tau'_1)=\{\emptyset,\{{\mathcal X}_{U_1},{\mathcal X}_{U_3}\},C(Y,{\bf S})\}$$ and $$t_{{\cal F}'_{1}}(\tau''_1)=\{\emptyset,\{{\mathcal X}_{U_2},{\mathcal X}_{U_3}\},C(Y,{\bf S})\}.$$

\bigskip

$\bullet$ Finally, let ${\cal F}_2=\{\{{\cal F}'_1\}\}$, where ${\cal F}'_1=\{\{U_1,U_3\},\{U_2,U_3\}\}$. Then, ${\cal O}^2({\cal F}_2)=\{\{{\cal F}_2\}\}$. We then consider the
topology $\tau_2=\{\emptyset,\{{\cal F}_2\}\}$ on ${\cal O}^2({\cal
F}_2)$ and define, using the {\bf Step 2}, the topology $t_{{\cal F}_{2}}(\tau_2)$ on the set $C(Y,{\bf S})$. We observe that $t_{{\cal F}_{2}}(\tau_2)=\{\emptyset,\{{\mathcal X}_{U_3}\},C(Y,{\bf S})\}$.

\bigskip
\noindent
{\bf (2)} Let ${\bf S}=\{0,1\}$ with the topology $\tau=\{\emptyset,\{1\},{\bf S}\}$ and let $Y=\{a,b\}$ with the topology $${\cal O}(Y)=\{\emptyset,U=\{a\},Y\}.$$  The Scott topology on ${\cal O}(Y)$ is the family $$\tau_{Sc}({\cal O}(Y))=\{\emptyset,\{Y\},\{U,Y\},{\cal O}(Y)\}.$$ Therefore, the Isbell topology on $C(Y,{\bf S})$ is the family $$t_{Is}=\{\emptyset,\{{\mathcal X}_Y\},\{{\mathcal X}_U,{\mathcal X}_Y\},C(Y,{\bf S})\}.$$ Let ${\cal F}_1={\cal P}({\cal O}(Y))=\{\emptyset,\{\emptyset\},\{U\},\{Y\},\{\emptyset,U\},\{\emptyset,Y\},\{U,Y\},{\cal O}(Y)\}$. Then, the set ${\cal O}^1({\cal F}_1)$ contains the following elements:

\smallskip

$O^1(\emptyset)=\{\{\emptyset\},\{\emptyset,U\},\{\emptyset,Y\},{\cal O}(Y)\}$,

\smallskip

$O^1(U)=\{\{U\},\{\emptyset,U\},\{U,Y\},{\cal O}(Y)\}$,

\smallskip

$O^1(Y)=\{\{Y\},\{\emptyset,Y\},\{U,Y\},{\cal O}(Y)\}$.

\medskip
\noindent
We consider the poset $({\cal O}^1({\cal F}_1),\subseteq)$. We observe that the Scott topology $\tau_{Sc}({\cal O}^1({\cal F}_1))$ is the discrete topology on ${\cal O}^1({\cal F}_1)$ and $t_{{\cal F}_{1}}(\tau_{Sc}({\cal O}^1({\cal F}_1)))$ is the discrete topology on $C(Y,{\bf S})$.

\bigskip
\noindent
{\bf (3)} Let ${\bf S}$ be the Sierpi\'{n}ski space, $Y=\mathbb R$ with the usual topology, 

\medskip

$D_1=\{(a,b)\subseteq\mathbb R:a,b\in\mathbb R,1\in(a,b)\}$,

$D_2=\{(a,b)\subseteq\mathbb R:a,b\in\mathbb R,2\in(a,b)\}$,

$\ldots$

$D_n=\{(a,b)\subseteq\mathbb R:a,b\in\mathbb R,n\in(a,b)\}$,

$\ldots$

\medskip
\noindent
and $${\cal F}_1=\{D_n:n=1,2,\ldots\}.$$ Then, $O^1(\emptyset)=O^1(\mathbb R)=\emptyset$. Moreover, for every $U\in {\cal O}(\mathbb R)$ we have $$O^1(U)=\begin{cases} \emptyset, & {\rm if} \ U \ {\rm is \ not \ bounded \ interval}\\
\emptyset, & {\rm if} \ \{1,2,\ldots\}\cap U=\emptyset \\ \{D_n:n\in\{1,2,\ldots\}\cap U\}, & {\rm otherwise}.\end{cases}$$ We consider the set ${\cal O}^1({\cal F}_1)=\{O^1(U):U\in {\cal O}(\mathbb R)\}$ and the
topology 
\begin{align*}
\tau_1=\,&\{\{O^1(0,2)\},\{O^1(0,2),O^1(2,4)\},\{O^1(0,2),O^1(2,4),O^1(4,6)\},\ldots\}\cup\\
&\{\emptyset\}\cup\{{\cal O}^1({\cal F}_1)\}
\end{align*}
on ${\cal O}^1({\cal F}_1)$. We observe that $$t_{{\cal F}_{1}}(\tau_1)=\{\emptyset,C(Y,{\bf S}),\{{\mathcal X}_{(0,2)}\},\{{\mathcal X}_{(0,2)},{\mathcal X}_{(2,4)}\},\{{\mathcal X}_{(0,2)},{\mathcal X}_{(2,4)},{\mathcal X}_{(4,6)}\},\ldots\}.$$

}
\end{example}

\section{Basic remarks for the ${\cal F}_n(\tau_n)$-family-open\\ topologies}

\begin{pro}{\rm Let $Y$ and $Z$ be two topological spaces,
$\tau_n$ an arbitrary  topology on ${\cal O}^n({\cal F}_n)$, where
$n=1,2,\ldots,$ and $\tau_{n+1}\subseteq{\cal P}({\cal
O}^{n+1}({\cal F}_{n+1}))$ the family, containing the empty set, which is defined
as follows $$I\!\! H\in\tau_{n+1} \ {\rm if \ and \ only \ if} \
\{O^{n}(U)\in{\cal O}^n({\cal F}_n):O^{n+1}(U)\in I\!\!
H\}\in\tau_{n}.$$ Then, the following statements are true:

(1) The family $\tau_{n+1}$ defines a topology on ${\cal
O}^{n+1}({\cal F}_{n+1})$.

(2) $t_{{\cal F}_{n+1}}(\tau_{n+1})\subseteq t_{{\cal F}_n}(\tau_n)$.

(3) If the topology $t_{{\cal F}_n}(\tau_n)$ on $C(Y,Z)$ is splitting, then
the topology $t_{{\cal F}_{n+1}}(\tau_{n+1})$ on $C(Y,Z)$ is splitting, too.

(4) If the topology $t_{{\cal F}_{n+1}}(\tau_{n+1})$ on $C(Y,Z)$ is
jointly continuous, then the topology $t_{{\cal F}_{n}(\tau_n)}$ on $C(Y,Z)$
is jointly continuous, too.}
\end{pro}
\noindent {\bf Proof.} (1) It is clear that $\emptyset, {\cal
O}^{n+1}({\cal F}_{n+1})\in \tau_{n+1}$. Let ${I\!\! H}_1,{I\!\!
H}_2\in\tau_{n+1}$. Then, $$I\!\! H'_1=\{O^{n}(U)\in{\cal O}^n({\cal
F}_n):O^{n+1}(U)\in I\!\! H_1\}\in\tau_n$$ and
$$I\!\! H'_2=\{O^{n}(U)\in{\cal O}^n({\cal F}_n):O^{n+1}(U)\in I\!\!
H_2\}\in\tau_{n}.$$ Since $\tau_n$ is a topology on ${\cal
O}^n({\cal F}_n)$ and $$I\!\! H'_1\cap I\!\! H'_2=\{O^{n}(U)\in{\cal
O}^n({\cal F}_n):O^{n+1}(U)\in I\!\! H_1\cap I\!\! H_2\},$$ $I\!\!
H'_1\cap I\!\! H'_2\in\tau_n$. Therefore, ${I\!\! H}_1\cap {I\!\!
H}_2\in\tau_{n+1}$.

\medskip
\noindent Now, let ${I\!\! H}_i\in\tau_{n+1}$, $i\in I$. Then,
$$I\!\! H'_i=\{O^{n}(U)\in{\cal O}^n({\cal F}_n):O^{n+1}(U)\in I\!\!
H_i\}\in\tau_{n}, \ i\in I.$$ Since $\tau_n$ is a topology on ${\cal
O}^n({\cal F}_n)$ and $$\bigcup_{i\in I}I\!\!
H'_i=\{O^{n}(U)\in{\cal O}^n({\cal F}_n):O^{n+1}(U)\in \bigcup_{i\in
I}I\!\! H_i\},$$ $\bigcup_{i\in I}I\!\! H'_i\in\tau_n$. Thus,
$\bigcup_{i\in I}I\!\! H_i\in\tau_{n+1}$. By the above, $\tau_{n+1}$
is a topology on ${\cal O}^{n+1}({\cal F}_{n+1})$. 

\medskip
\noindent (2) Let $g\in C(Y,Z)$ and $<I\!\! H, U>_{n+1}$ be a
subbasic open set of the topology $t_{{\cal F}_{n+1}}(\tau_{n+1})$, such that
$g\in <I\!\! H, U>_{n+1}$. Consider the set $$I\!\!
H'=\{O^{n}(U)\in{\cal O}^n({\cal F}_n):O^{n+1}(U)\in I\!\! H\}.$$
Then, $I\!\! H'\in\tau_{n}$. We prove that $$g\in <I\!\! H',
U>_{n}\subseteq <I\!\! H, U>_{n+1}.$$ First, we prove that $g\in
<I\!\! H', U>_{n}$. Since $$g\in <I\!\! H, U>_{n+1}=\{f\in
C(Y,Z):O^{n+1}(f^{-1}(U))\in I\!\! H\},$$ we have that $O^{n+1}(g^{-1}(U))\in
I\!\! H$ and, therefore, $O^{n}(g^{-1}(U))\in I\!\! H'$. Hence,
$$g\in <I\!\! H', U>_{n}=\{f\in C(Y,Z):O^{n}(f^{-1}(U))\in I\!\!
H'\}.$$ Now, we prove that $<I\!\! H', U>_{n}\subseteq <I\!\! H,
U>_{n+1}$. Let $\varphi\in <I\!\! H', U>_{n}$. Then,
$O^{n}(\varphi^{-1}(U))\in I\!\! H'$ and, hence,
$O^{n+1}(\varphi^{-1}(U))\in I\!\! H$. Thus, $\varphi\in <I\!\! H,
U>_{n+1}$.

\medskip
\noindent The proofs for (3) and (4) are a straightforward
verification of statement (2). $\Box$

\begin{pro}{\rm
Let $\tau_n$ be a topology on the set ${\cal O}^n({\cal F}_n)$. Then, the mapping $h_n:C_{t_{{\cal F}_n}(\tau_n)}(Y,{\bf S})\to({\cal O}^n({\cal F}_n),\tau_n)$, defined by $h_n(f)=O^n(f^{-1}(\{1\}))$, for every $f\in C_{t_{{\cal F}_n}(\tau_n)}(Y,{\bf S})$, is continuous, open and onto.}
\end{pro}
{\bf Proof.} First we prove that the mapping $h_n$ is continuous. Let
$\H\in\tau_n$. Then, $$h_n^{-1}(\H)=\{f\in C_{t_{{\cal F}_n}(\tau_n)}(Y,{\bf
S}):O^n(f^{-1}(\{1\}))\in\H\}=<\H,\{1\}>_{n},$$ which shows continuity. Now, we prove that
the mapping $h_n$ is open. Let $<\H,\{1\}>_{n}$ be a subbasic open set
in $C_{t_{{\cal F}_n}(\tau_n)}(Y,{\bf S})$. It suffices to prove that
$h_n(<\H,\{1\}>_{n})=\H$. Indeed, let $O^n(U)\in
h_n(<\H,\{1\}>_{n})$. Then, there exists $f\in <\H,\{1\}>_{n}$, such
that $O^n(U)=h_n(f)=O^n(f^{-1}(\{1\}))\in\H$. Hence,
$h_n(<\H,\{1\}>_{n})\subseteq\H$. Let $O^n(U)\in\H$. Consider the
characteristic function ${\mathcal X}_U:Y\to{\bf S}$, of $U$, defined
by $${\mathcal X}_U(y)=\begin{cases}1 & {\rm if} \ y\in U,\\ 0 &
{\rm if} \ y\in Y\setminus U.\end{cases}$$ Then, ${\mathcal
X}_U^{-1}(\{1\})=U$ and, therefore, $$h_n({\mathcal
X}_U)=O^n({\mathcal X}_U^{-1}(\{1\}))=O^n(U)\in\H.$$ It follows that
$O^n(U)\in h_n(<\H,\{1\}>_{n})$. Hence, $\H\subseteq
h_n(<\H,\{1\}>_{n})$. So, we have that $h_n(<\H,\{1\}>_{n})=\H$.
Similarly, we can prove that $h_n$ is onto. $\Box$

\begin{notation}
{\rm By $\Phi_{{\cal F}_{n}}$ we denote the map from ${\cal
O}(Y)$ to ${\cal O}^{n}({\cal F}_{n})$, for which $\Phi_{{\cal
F}_{n}}(U)=O^{n}(U)$, for every $U\in{\cal O}(Y)$.}
\end{notation}

\begin{pro}{\rm
Let $\tau_n$ be a topology on the set ${\cal O}^n({\cal F}_n)$. If
the map $\Phi_{{\cal F}_{n}}$ is one-to-one, then the mapping
$h_n:C_{t_{{\cal F}_n}(\tau_n)}(Y,{\bf S})\to({\cal O}^n({\cal
F}_n),\tau_n)$ is a homeomorphism.}
\end{pro}
{\bf Proof.} Suppose that the map $\Phi_{{\cal F}_{n}}$ is
one-to-one. By Proposition 3.2 it suffices to prove that the map
$h_n$ is one-to-one. Indeed, let $f,g\in C_{t_{{\cal F}_n}(\tau_n)}(Y,{\bf
S})$ such that $h_n(f)=h_n(g)$. Then,
$O^n(f^{-1}(\{1\}))=O^n(g^{-1}(\{1\}))$ or, equivalently, $\Phi_{{\cal
F}_{n}}(f^{-1}(\{1\}))=\Phi_{{\cal F}_{n}}(g^{-1}(\{1\}))$. By
assumption, $f^{-1}(\{1\})=g^{-1}(\{1\})$. Therefore, $f=g$. $\Box$

\begin{pro}{\rm Let $Y$ be a topological space, $\tau_0$ a ${\rm T}_0$-topology on ${\cal O}(Y)$, and $\tau_i$ a ${\rm T}_0$-topology on the set ${\cal O}^i({\cal F}_i)$, $i=1,\ldots,n-1$. If ${\cal F}_{i+1}=\tau_i$ for every $i=0,1,\ldots,n-1$, then the spaces $C_{t_{{\cal F}_n}(\tau_n)}(Y,{\bf S})$ and $({\cal O}^n({\cal F}_n),\tau_n)$ are homeomorphic to each other.}
\end{pro}
{\bf Proof.} By Proposition 3.3 it suffices to prove that the map
$\Phi_{{\cal F}_{n}}$ is one-to-one. Let $U,V\in {\cal O}(Y)$ such
that $U\neq V$. We prove that $\Phi_{{\cal F}_{n}}(U)\neq
\Phi_{{\cal F}_{n}}(V)$. Indeed, following our construction we have
the following steps:

\medskip
\noindent {\bf Step 1.} Since ${\cal F}_{1}=\tau_0$, the families
$O^1(U)$ and $O^1(V)$ are the filters of open neighborhoods of $U$
and $V$, respectively, in the topology $\tau_0$. Since the topology
$\tau_0$ on the set ${\cal O}(Y)$ is ${\rm T}_0$, we get that
$O^1(U)\neq O^1(V)$.

\medskip
\noindent {\bf Step 2.} Since ${\cal F}_{2}=\tau_1$, the families
$O^2(U)$ and $O^2(V)$ are the filters of open neighborhoods of $O^1(U)$ and $O^1(V)$, respectively, in the topology $\tau_1$. Since the topology
$\tau_1$ on the set ${\cal O}^1({\cal F}_1)$ is ${\rm T}_0$,
$O^2(U)\neq O^2(V)$.

\medskip
\noindent We continue in the same manner to a Step 3, Step 4, etc.
Step $n$ will be as follows.

\medskip
\noindent {\bf Step n.} Since ${\cal F}_{n}=\tau_{n-1}$, the
families $O^{n}(U)$ and $O^{n}(V)$ are the filters of open
neighborhoods of $O^{n-1}(U)$ and $O^{n-1}(V)$, respectively, in the topology
$\tau_{n-1}$. Since the topology $\tau_{n-1}$ on the set ${\cal
O}^{n-1}({\cal F}_{n-1})$ is ${\rm T}_0$, $O^n(U)\neq O^n(V)$ or,
equivalently, $\Phi_{{\cal F}_{n}}(U)\neq \Phi_{{\cal F}_{n}}(V)$. $\Box$

\section{A characterization of ${\cal A}$-splitting topologies for  the ${\cal F}_n(\tau_n)$-family-open topologies on\\ $C(Y,Z)$}

\begin{definition}
{\rm Let $\tau_n$ be a topology on the set ${\cal O}^n({\cal F}_n)$.
We say that a map $M:X\times {\cal O}(Z)\to{\cal
O}^n({\cal F}_n)$ is {\it continuous with respect to the first
variable} if, for every fixed element $U$ of ${\cal O}(Z)$, the map
$M_{U}: X\to({\cal O}^n({\cal F}_n),\tau_n)$, for which
$M_{U}(x)=M(x,U)$ for every $x\in X$, is
continuous.}
\end{definition}

\begin{notation}
{\rm Let $X$ be an arbitrary topological space.  If $F:X\times Y\to
Z$ is a continuous map, then by $\overline{F}_n$ we denote the map
from $X\times {\cal O}(Z)$ to ${\cal O}^n({\cal F}_n)$ for which
$$\overline{F}_n(x,U)=O^n(F_x^{-1}(U)),$$ for every $x\in X$ and
$U\in {\cal O}(Z)$.}
\end{notation}

\begin{pro}
{\rm Let ${\cal A}$ be an arbitrary family of topological spaces and
$\tau_n$ a topology on ${\cal O}^n({\cal F}_n)$, where $n=1,2,\ldots.$
The topology $t _{{\cal F}_n(\tau_n)} $ on $C(Y,Z)$ is ${\cal
A}$-splitting if and only if, for every space $X\in {\cal A}$, the
continuity of the map $F:X\times Y\to Z$ implies the continuity of
the map $\overline{F}_n:X\times {\cal O}(Z)\to {\cal O}^n({\cal
F}_n)$, with respect to the first variable.}
\end{pro}
\noindent {\bf Proof.} Suppose that the topology $t_{{\cal F}_n}(\tau_n)$ on
$C(Y,Z)$ is ${\cal A}$-splitting, $X$ is an element of ${\cal A}$
and $F:X\times Y\to Z$ is a continuous map. We must prove that the
map $\overline{F}_n:X\times {\cal O}(Z)\to {\cal O}^n({\cal F}_n)$
is continuous with respect to the first variable. For this, let $U$ be a fixed element of ${\cal O}(Z)$, $x\in X$,
$I\!\! H\in \tau_n$, and also let
$$(\overline{F}_{n})_U(x)=\overline{F}_n(x,U)=O^n(F^{-1}_x(U))\in
I\!\!H.$$ We need to find an open neighborhood $V$ of $x$, in $X$,
such that $$(\overline{F}_{n})_U(V)=\overline{F}_n(V\times\{U\})\s
I\!\!H.$$ We consider the open set $<I\!\! H, U>_n$ of the space $C_{t
_{{\cal F}_n}(\tau_n)}(Y,Z)$. Then $$F_x\in <I\!\! H, U>_n.$$ Since the
topology $t_{{\cal F}_n}(\tau_n)$ is ${\cal A}$-splitting, the map
$\h{F}:X\to C_{t_{{\cal F}_n}(\tau_n)}(Y,Z)$ is continuous. Thus,
there exists an open neighborhood $V$ of $x$, in $X$, such that
$\h{F}(V)\s <I\!\! H,U>_n$. We will prove that
$$(\ov{F}_{n})_U(V)=\ov{F}_n(V\times\{U\})\s I\!\! H.$$ Indeed, let
$x^{\prime}\in V$. Then $$\h{F}(x^{\prime})=F_{x^{\prime}}\in
<I\!\! H, U>_n.$$ But this implies that $O^n(F^{-1}_{x^{\prime}}(U))\in
I\!\! H$ or, equivalently, that $$(\ov{F}_{n})_U(x')=\ov{F}_n(x^{\prime},
U)\in I\!\! H.$$ Conversely, suppose that for every element $X$, of
${\cal A}$, the continuity of the map $F:X\times Y\to Z$ implies the
continuity of the map $\overline{F}_n$, with respect to the first
variable. We will prove that the topology $t_{{\cal F}_n}(\tau_n)$ on $C(Y,Z)$ is
${\cal A}$-splitting. For this, let $X\in\mathcal A$ and let  also $F:X\times Y\to Z$ be a
continuous map. We need to prove that the map $\h{F}:X\to
C_{t_{{\cal F}_n}(\tau_n)} (Y,Z)$ is continuous. Indeed, let $x\in X$
and let $<I\!\! H, U>_n$ be an open neighborhood of $\h{F}(x)$, in
$C_{t_{{\cal F}_n}(\tau_n)}(Y,Z)$. It suffices to prove that there
exists an open neighborhood $V$ of $x$, in $X$, such that
$\h{F}(V)\subseteq <I\!\! H, U>_n$. But, since $\h{F}(x)=F_x\in <I\!\! H,
U>_n$, we have that $O^n((\h{F}(x))^{-1}(U))\in I\!\! H$, that is
$$(\ov{F}_{n})_U(x)=\ov{F}_n(x, U)\in I\!\! H.$$ By assumption, the
map $\ov{F}_n$ is continuous with respect to first variable. Thus,
there exists an open neighborhood $V$ of $x$, in $X$, such that
$$(\ov{F}_{n})_U(V)=\ov{F}_n(V\times\{U\})\s I\!\! H.$$ We will now prove
that $\h{F}(V)\s <I\!\! H, U>_n$. Indeed, let $x^{\prime}\in V$.
Then $$(\ov{F}_{n})_U(x^{\prime})=\ov{F}_n(x^{\prime},U)\in I\!\!
H,$$ that is, $O^n(F^{-1}_{x^{\prime}}(U))\in I\!\! H$ and,
therefore, $\h{F}(x^{\prime})\in <I\!\! H, U>_n$. $\Box$

\begin{co}
{\rm Let $\tau_n$ a topology on ${\cal O}^n({\cal F}_n)$, where
$n=1,2,\ldots.$ Then, the topology $t_{{\cal F}_n}(\tau_n)$ on $C(Y,Z)$ is
splitting, if and only if for every space $X$ the continuity of the
map $F:X\times Y\to Z$ implies the continuity of the map
$\overline{F}_n:X\times {\cal O}(Z)\to {\cal O}^n({\cal F}_n)$, with
respect to the first variable.}
\end{co}

\begin{notation}
{\rm For every continuous map $F:X\times Y\to Z$ we denote by $F^{*}$ the
map from $X\times{\cal O}(Z)$ to ${\cal O}(Y)$, for which
$$F^{*}(x,U)=F_x^{-1}(U),$$ for every $(x,U)\in X\times{\cal O}(Z)$. }
\end{notation}

\begin{definition}
{\rm Let $F:X\times Y\to Z$ be a continuous map and $\tau_{Sc}({\cal O}(Y))$ be the
Scott topology on ${\cal O}(Y)$. We say that a map
$F^{*}:X\times{\cal O}(Z)\to {\cal O}(Y)$ is {\it Scott continuous with
respect to the first variable} if, for every fixed element $U$ of
${\cal O}(Z)$, the map $(F^{*})_{U}: X\to {\cal O}(Y)$, for which
$(F^{*})_{U}(x)=F_x^{-1}(U)$ for every $x\in X$, is continuous.}
\end{definition}

\begin{remark}
{\rm We observe that for every $n=1,2,\ldots$ the pair $({\cal
O}^{n}({\cal F}_{n}),\subseteq)$ is a partially ordered set. This allows us to consider the Scott topology on the set ${\cal O}^{n}({\cal
F}_{n})$.}
\end{remark}

\begin{definition}
{\rm We say that a map $f:{\cal O}(Y)\to {\cal O}^{n}({\cal F}_{n})$
is {\it Scott continuous}, if the map $f$ is continuous when the sets ${\cal
O}(Y)$ and ${\cal O}^{n}({\cal F}_{n})$ are endowed with the Scott topology.}
\end{definition}

\begin{pro}
{\rm Let ${\cal A}$ be an arbitrary family of topological spaces. If
the map $\Phi_{{\cal F}_{n}}:{\cal O}(Y)\to{\cal O}^{n}({\cal
F}_{n})$ is Scott continuous, then the topology $t_{{\cal F}_{n}}(\tau_{Sc}({\cal O}^{n}({\cal
F}_{n})))$
on $C(Y,Z)$ is ${\cal A}$-splitting, for every topological space
$Z$.}
\end{pro}
\noindent {\bf Proof.} Let $X\in{\cal A}$ and let $F:X\times Y\to Z$ be
a continuous map. Using Proposition 4.1 it suffices to prove that the map
$\ov{F}_{n}:X\times{\cal O}(Z)\to{\cal O}^{n}({\cal F}_{n})$ is
continuous with respect to the first variable. Indeed, let
$U\in{\cal O}(Z)$. For every $x\in X$ we have
$$\ov{F}(x,U)=O^{n}(F_{x}^{-1}(U))=(\Phi_{{\cal F}_{n}}\circ F^{*})(x,U).$$
So, it suffices to prove that the map $F^{*}$ is Scott continuous with
respect to the first variable. Let $x\in X$, let $I\!\! H$ be a Scott open set in
${\cal O}(Y)$ and let $(F^{*})_{U}(x)\in I\!\! H$. We will find an
open neighborhood $V$ of $x\in X$, such that $(F^{*})_{U}(V)\subseteq
I\!\! H$. But
$$(F^{*})_{U}(x)=F^{-1}_{x}(U)=(\h{F}(x))^{-1}(U)\in I\!\! H,$$ that is
$\h{F}(x)\in (I\!\! H, U)$. Since the Isbell topology on the set
$C(Y,Z)$ is splitting, the map $\h{F}:X\to C_{t_{Is}}(Y,Z)$ is continuous.
Hence, there exists an open neighborhood $V$ of $x\in X$, such that
$\h{F}(V)\subseteq (I\!\! H, U)$. We now prove that $(F^{*})_{U}(V)\subseteq
I\!\! H$. Indeed, let $y\in V$. Then $$\h{F}(y)\in (I\!\! H, U) \ {\rm or} \ (\h{F}(y))^{-1}(U)\in I\!\! H$$ which is equivalent to say that $(F^{*})_U(y)\in I\!\! H$. Therefore, the map $F^{*}$ is Scott
continuous with respect to the first variable and this finishes the proof. $\Box$

\begin{co}
{\rm If the map $\Phi_{{\cal F}_{n}}:{\cal O}(Y)\to{\cal
O}^{n}({\cal F}_{n})$ is Scott continuous, then the topology
$t_{{\cal F}_{n}}(\tau_{Sc}({\cal O}^{n}({\cal
F}_{n})))$ on $C(Y,Z)$ is splitting, for any topological
space $Z$.}
\end{co}

\begin{lemma}
{\rm Let $\tau_{Sc}({\cal O}^{n}({\cal F}_{n}))$ be the
Scott topology on ${\cal O}^{n}({\cal F}_{n})$ and let ${\cal
F}_{n+1}\subseteq \tau_{Sc}({\cal O}^{n}({\cal F}_{n}))$. If $O^{n}(V)\subseteq O^{n}(U)$, then $O^{n+1}(V)\subseteq O^{n+1}(U)$.}
\end{lemma}
\noindent
{\bf Proof.} Let $$\psi\in O^{n+1}(V)=\{\varphi\in{\cal F}_{n+1}:O^{n}(V)\in\varphi\}.$$ Then, $O^{n}(V)\in\psi$,
$O^{n}(V)\subseteq O^{n}(U)$ and $\psi$ is a Scott open set in ${\cal O}^{n}({\cal F}_{n})$. Therefore, $O^{n}(U)\in\psi$.
This means that $$\psi\in O^{n+1}(U)=\{\varphi\in{\cal F}_{n+1}:O^{n}(U)\in\varphi\},$$ which finishes the proof. $\Box$

\begin{lemma}
{\rm Let $\tau_{Sc}({\cal O}(Y))$ be the
Scott topology on ${\cal O}(Y)$ and let ${\cal
F}_{1}\subseteq \tau_{Sc}({\cal O}(Y))$. If $\{V_{\alpha}:\alpha\in A\}\subseteq {\cal O}(Y)$, then $$O^{n}(\bigcup_{\alpha\in A}V_{\alpha})=\bigcup_{\lambda\in\Lambda}O^{n}(\bigcup_{\alpha\in\lambda}V_{\alpha}),$$ where $\Lambda$ is the set of all finite subsets of $A$.}
\end{lemma}
\noindent
{\bf Proof.} By Lemma 4.1 it is clear that $$\bigcup_{\lambda\in\Lambda}O^{1}(\bigcup_{\alpha\in\lambda}V_{\alpha})\subseteq O^{1}(\bigcup_{\alpha\in A}V_{\alpha}).$$ We prove that $$O^{1}(\bigcup_{\alpha\in
A}V_{\alpha})\subseteq\bigcup_{\lambda\in\Lambda}O^{1}(\bigcup_{\alpha\in\lambda}V_{\alpha}).$$ Indeed, let
$$\psi\in O^{1}(\bigcup_{\alpha\in A}V_{\alpha})=\{\varphi\in{\cal F}_{1}:\bigcup_{\alpha\in A}V_{\alpha}\in\varphi\}.$$ Then, $\bigcup_{\alpha\in A}V_{\alpha}\in\psi$.
Since $\psi$ is a Scott open set in ${\cal O}(Y)$, there exists an
element $\lambda_0\in\Lambda$ such that $\bigcup_{\alpha\in\lambda_0}V_{\alpha}\in\psi$. Hence $$\psi\in O^{1}(\bigcup_{\alpha\in\lambda_0}V_{\alpha})=\{\varphi\in{\cal F}_{1}:\bigcup_{\alpha\in \lambda_0}V_{\alpha}\in\varphi\}\subseteq\bigcup_{\lambda\in\Lambda}O^{1}(\bigcup_{\alpha\in\lambda}V_{\alpha}). \ \Box$$

\begin{pro}
{\rm Let ${\cal A}$ be an arbitrary family of spaces, $\tau_{Sc}({\cal O}(Y))$ the
Scott topology on ${\cal O}(Y)$ and let ${\cal
F}_{1}\subseteq \tau_{Sc}({\cal O}(Y))$. Then, the topology
$t_{{\cal F}_{1}}(\tau_{Sc}({\cal O}^{1}({\cal
F}_{1})))$ on $C(Y,Z)$ is ${\cal A}$-splitting.}
\end{pro}
\noindent {\bf Proof.} By Proposition 4.2, it suffices to prove that the
map $$\Phi_{{\cal F}_{1}}:{\cal O}(Y)\to {\cal O}^{1}({\cal
F}_{1})$$ is Scott continuous. Let $I\!\! H$ be a Scott open set on
${\cal O}^{1}({\cal F}_{1})$. We prove that the set
$$I\!\! H'=\Phi^{-1}_{{\cal F}_{1}}(I\!\! H)=\{V\in{\cal O}(Y):\Phi_{{\cal F}_{1}}(V)=O^{1}(V)\in I\!\! H\}$$
is Scott open in ${\cal O}(Y)$. For this, it suffices to prove that $I\!\! H'$ satisfies the two conditions from the definition of Scott topology.

\medskip
\noindent (1) Let $V\in I\!\! H',\ U\in{\cal O}(Y)$ and $V\subseteq
U$. We prove that $U\in I\!\! H'$. Since $V\in I\!\! H'$, $O^{1}(V)\in I\!\! H$. Also, since $V\subseteq U$, by Lemma 4.1, $O^{1}(V)\subseteq O^{1}(U)$. Thus,
$O^{1}(U)\in I\!\! H$ and, therefore, $U\in I\!\! H'$.

\medskip
\noindent (2) Now, let $V=\bigcup_{\alpha\in A}V_{\alpha}\in I\!\! H'$, where $V_{\alpha}\in{\cal O}(Y),\ \alpha\in A$. Then $$O^{1}(V)=O^{1}(\bigcup_{\alpha\in A}V_{\alpha})\in I\!\! H.$$ By Lemma 4.2,
$$O^{1}(\bigcup_{\alpha\in
A}V_{\alpha})=\bigcup_{\lambda\in\Lambda}O^{1}(\bigcup_{\alpha\in\lambda}V_{\alpha}),$$
where $\Lambda$ is the set of all finite subsets of $A$.
Since $I \!\! H$ is a Scott open set, there exist
$\lambda_1,\ldots,\lambda_n$ elements of $\Lambda$ such that
$$\bigcup_{i=1}^{n}O^1(\bigcup_{\alpha\in\lambda_{i}}V_{\alpha})\in I\!\! H.$$
Also, since
$$\bigcup_{i=1}^{n}O^{1}(\bigcup_{\alpha\in\lambda_{i}}V_{\alpha})\subseteq
O^1(\bigcup_{i=1}^{n}(\bigcup_{\alpha\in\lambda_{i}}V_{\alpha}))$$
and $I\!\!H$ is a Scott open set in $O^{1}({\cal F}_{1})$, we have
that
$$O^{1}(\bigcup_{i=1}^{n}(\bigcup_{\alpha\in\lambda_{i}}V_{\alpha}))\in I\!\!H,$$
which means that
$$\bigcup_{i=1}^{n}(\bigcup_{\alpha\in\lambda_{i}}V_{\alpha})\in I\!\!H'.$$
Therefore, $I\!\!H'=\Phi_{{\cal F}_{1}}^{-1}(I\!\!H)$ is a Scott open set. $\Box$

\begin{co}
{\rm Let ${\cal F}_{1}\subseteq \tau_{Sc}({\cal O}(Y))$. Then, the
topology $t_{{\cal F}_{1}}(\tau_{Sc}({\cal O}^{1}({\cal
F}_{1})))$ on $C(Y,Z)$ is splitting.}
\end{co}

\section{A characterization of
${\cal A}$-jointly continuous topologies for  the ${\cal
F}_n(\tau_n)$-family-open topologies on $C(Y,Z)$}

\begin{notation}
{\rm Let $X$ be an arbitrary topological space. If $G:X\to C(Y,Z)$
is a map, then by $\overline{G}_n$, where $n=1,2,\ldots$, we denote
the map from $X\times {\cal O}(Z)$ to ${\cal O}^n({\cal F}_n)$, for
which $$\overline{G}_n(x,U)=O^n((G(x))^{-1}(U)),$$ for every $x\in
X$ and $U\in {\cal O}(Z)$.}
\end{notation}

\begin{pro}
{\rm Let ${\cal A}$ be an arbitrary family of topological spaces and
$\tau_n$ a topology on ${\cal O}^n({\cal F}_n)$, where $n=1,2,\ldots.$
The topology $t_{{\cal F}_n}(\tau_n)$, on $C(Y,Z)$, is ${\cal A}$-jointly
continuous, if and only if for every space $X\in {\cal A}$, the
continuity of the map $\overline{G}_n$, with respect to the first
variable, implies the continuity of the map $\t{G}:X\times Y\to Z$.}
\end{pro}
\noindent
{\bf Proof.} Suppose that the topology $t_{{\cal F}_n}(\tau_n)$ on
$C(Y,Z)$ is ${\cal A}$-jointly continuous, $X$ is an element of
${\cal A}$ and let $\overline{G}_n:X\times {\cal O}(Z)\to {\cal
O}^n({\cal F})$ be a continuous map, with respect to the first
variable. We will prove that the map $\t{G}:X\times Y\to Z$ is
continuous. But, since the topology $t_{{\cal F}_n}(\tau_n)$ on $C(Y,Z)$ is ${\cal
A}$-jointly continuous, it  suffices to prove that the map $G:X\to
C_{t_{{\cal F}_n}(\tau_n)}(Y,Z)$ is continuous. Indeed, let $x\in X$, $U\in
{\cal O}(Z)$ and let $I\!\! H\in \tau_n$, such that $G(x)\in <I\!\!
H,U>_n$. We need to find an open neighborhood $V$ of $x$, in $X$,
such that $G(V)\s <I\!\! H,U>_n$. But it is true that
$$O^n((G(x))^{-1}(U))\in I\!\! H$$ and, since the map
$\overline{G}_n:X\times {\cal O}(Z)\to {\cal O}^n({\cal F})$ is
continuous with respect to the first variable and
$O^n((G(x))^{-1}(U))\in I\!\! H$, there exists an open
neighborhood $V$ of $x$, in $X$, such that
$$(\overline{G}_{n})_U(V)=\overline{G}_n(V\times\{U\})\s I\!\! H.$$ It
remains to prove that $G(V)\s <I\!\! H,U>_n$. Indeed, let
$x^{\prime}\in V$. Then $$\ov{G}_n(x^{\prime},U)=
O^n((G(x^{\prime}))^{-1}(U))\in I\!\! H$$ or, equivalently,
$G(x^{\prime})\in < I\!\! H,U>_n$.

\medskip
\noindent Conversely, suppose that for every space $X\in {\cal A}$
the continuity of the map $\overline{G}_n$, with respect to the first
variable, implies the continuity of the map $$\t{G}:X\times Y\to
Z.$$ We prove that $t_{{\cal F} _n}(\tau_n)$ is ${\cal A}$-jointly
continuous. Let $X\in {\cal A}$ and let $$G:X\to C_{t_{{\cal F} _n}(\tau_n)}
(Y,Z)$$ be a continuous map. We prove that the map $\t{G}:
X\times Y \to Z$ is continuous. For this, it suffices to prove that
the map $$\ov{G}_n:X\times {\cal O}(Z)\to {\cal O}^n({\cal F})$$ is
continuous with respect to the first variable. Indeed, let $x\in X$,
$U\in {\cal O}(Z)$ and let $I\!\! H\in\tau_n$, so that
$$(\ov{G}_{n})_U(x)=\ov{G}_n(x,U)\in I\!\! H.$$ We must find an open
neighborhood $V$ of $x$, in $X$, such that $\ov{G}_{n,U}(V)\s I\!\!
H$. By considering the open set $<I\!\! H,U>_n$, of the space
$C_{t_{{\cal F} _n}(\tau_n)}(Y,Z)$, we get that $G(x)\in <I\!\! H,U>_n$.
But, since the map $G:X\to C_{t_{{\cal F}_n}(\tau_n)}(Y,Z)$ is continuous
and $G(x)\in <I\!\! H,U>_n$, there exists an open neighborhood
$V$ of $x$, in $X$, such that $G(V)\s <I\!\! H,U>_n$. We now prove
that $$(\ov{G}_{n})_U(V)=\ov{G}_n(V\times\{U\})\s I\!\! H.$$ Indeed,
let $x^{\prime}\in V$. Then, $G(x^{\prime})\in <I\!\! H,U>_n$ or,
equivalently $$O^n((G(x^{\prime}))^{-1}(U))\in I\!\! H.$$ Thus,
$\ov{G}_n(x^{\prime},U)\in I\!\! H$ and, therefore
$\ov{G}_n(V\times\{U\})\s I\!\! H$. $\Box$

\begin{co}
{\rm Let $\tau_n$ be a topology on ${\cal O}^n({\cal F}_n)$, where
$n=1,2,\ldots.$ The topology $t_{{\cal F}_n}(\tau_n) $ on $C(Y,Z)$ is
jointly continuous if and only if, for every space $X$, the
continuity of the map $\overline{G}_n$, with respect to the first
variable, implies the continuity of the map $\t{G}:X\times Y\to Z$.}
\end{co}

\section{Some Open Problems}

In the past years, there has been a great deal of progress in the field of function spaces. There are several papers on this area (see, for example, \cite{DM}- \cite{DGL}, \cite{E1}, \cite{E2}, \cite{GEORGIOU}, \cite{K}, \cite{N}, \cite{XY}). In this section we give some problems concerning ${\cal F}_n(\tau_n)$-family-open topologies.

\medskip
\noindent
1. Under what conditions on the spaces $Y$ and $Z$ is the topology
$t_{{\cal F}_n}(\tau_n)$ splitting?

\medskip
\noindent
2. Under what conditions on the spaces $Y$ and $Z$ is the topology
$t_{{\cal F}_n}(\tau_n)$ jointly continuous?

\medskip
\noindent
3. Is the greatest splitting topology, which always exists, a ${\cal F}_n(\tau_n)$-family-open topology?

\medskip
\noindent
4. Let $t$ be an arbitrary set-open topology on $C(Y,Z)$.
Is this topology a ${\cal F}_n(\tau_n)$-family-open topology?

\medskip
\noindent
5. Find Arzela-Ascoli theorems  for the topology
$t_{{{\cal F}}_n}$, $n=1,2,\dots$.

\medskip
\noindent
6. Let $t_{{\cal F}_n}(\tau_n), n=1,2,\dots$ be the topologies on
$C(Y,Z)$ which we constructed on $C(Y,Z)$, in Proposition 3.1. Does there exist a
positive number $n$ such that: $$t_{{\cal F}_n}(\tau_n)=t_{{\cal F}_{n+1}}(\tau_{n+1})=t_{{\cal F}_{n+2}}(\tau_{n+2})=\dots \ ?$$

\medskip
\noindent 
7. Let ${\bf S}$ be the Sierpi\'{n}ski space. Then, the
set $C(Y,{\bf S})$ coincides with the set ${\mathcal CL}(Y)$, of
closed subsets of $Y$. Consider, on $C(Y,{\bf S})$, the
following topologies:\\
(a) the Vietoris topology $t_V$ and\\
(b) the Fell topology $t_F$.\\
Are the above topologies ${\cal F}_n(\tau_n)$-family-open topologies?

\medskip
\noindent 
8. Consider, on $C(Y,Z)$, the
following topologies (see, for example, \cite{DM}):\\
(a) the fine topology,\\
(b) the graph topology, and\\
(c) the Krikorian topology.\\
Are the above topologies ${\cal F}_n(\tau_n)$-family-open topologies?

\medskip
\noindent 
9. Let $X$ be an arbitrary space. By $w(X)$ denote the
weight of $X$. Is the following relation valid $w(C_{t_{{\cal F}_n}(\tau_n)}(Y,Z))\leq w(Y)w(Z)$, for $n=1,2,\ldots$?

\bigskip
\noindent
{\bf Acknowledgements.} The authors would like to thank the anonymous reviewers for their valuable comments and suggestions to improve the quality of the paper.

\bigskip

Email addresses:\\
georgiou@math.upatras.gr (Dimitris Georgiou)\\
thanasismeg13@gmail.com (Athanasios Megaritis)\\
kyriakos.papadopoulos1981@gmail.com (Kyriakos Papadopoulos)\\
bpetrop@master.math.upatras.gr (Vasilios Petropoulos)


\begin{thebibliography}{99}
\bibitem{ARENS} R. Arens, {\sl A topology for spaces of transformations},
Ann. of Math. 47(1946), 480--495.

\bibitem{AREDUG} R. Arens and J. Dugundji, {\sl Topologies for function
spaces}, Pacific J. Math. 1(1951),  5--31.

\bibitem{DM} G. Di Maio, L. Hol\'{a}, D. Hol\'{y} and R. McCoy, {\sl Topologies on the set space of continuous functions}, Topology Appl. 86 (1998), no. 2, 105--122.

\bibitem{DMN1} G. Di Maio, E. Meccariello, S. Naimpally, {\sl Hyper-continuous convergence in function spaces}, Questions and Answers in General Topology, vol. 22, no. 2, pp. 157–-162, 2004

\bibitem{DMN2} G. Di Maio, E. Meccariello, S. Naimpally, {\sl Hyper-continuous convergence in function spaces}, II, Ricerche di Matematica, vol. 54, no. 1, pp. 245–-254, 2005.

\bibitem{DMN3} G. Di Maio, E. Meccariello, S. Naimpally, {\sl Duality in function spaces}, Mediterr. J. Math. 3 (2006), no. 2, 189–-204.

\bibitem{DMN4} G. Di Maio, E. Meccariello, S. Naimpally, {\sl Hyperspace and function space are duals}, Questions Answers Gen. Topology 25 (2007), no. 1, 23–-43.

\bibitem{DGL} S. Dolecki, G. H. Greco, A. Lechicki, {\sl When do the upper Kuratowski topology (homeomorphically, Scott topology) and the co-compact topology coincide?} Trans. Amer. Math. Soc. 347 (1995), no. 8, 2869-–2884.

\bibitem{DUG} J. Dugundji, {\sl Topology}, Allyn and Bacon, Boston, Mass., 1966.

\bibitem{E1} M. H. Escardo, {\sl Function-space compactifications of function spaces}, Topology Appl. 120 (2002), no. 3, 441–-463.

\bibitem{E2} M. Escardo, R. Heckmann {\sl Topologies on spaces of continuous functions}, Proceedings of the 16th Summer Conference on General Topology and its Applications (New York), Topology Proc. 26 (2001/02), no. 2, 545–-564.

\bibitem{FOX} R. H. Fox, {\sl On topologies for function spaces}, Bull. Amer. Math. Soc.
51(1945), 429--432.

\bibitem{GEORGIOU} D. N. Georgiou,
S. D. Iliadis, and F. Mynard, {\sl Function Space Topologies},
(appear in Open Problems in Topology 2 (Elsevier)).

\bibitem{GEORGIOU1} D. N. Georgiou, S. D. Iliadis, and B. K. Papadopoulos,  {\sl Topologies on function spaces}, Studies in Topology, VII, Zap. Nauchn. Sem.
S.-Peterburg Otdel. Mat. Inst. Steklov (POMI), 208 (1992),  82-97.
J. Math. Sci., New York 81, (1996), No. 2, pp. 2506--2514.

\bibitem{SCOT} G. Gierz, K. H. Hofmann, K. Keimel, J. D. Lawson, M.
Mislove, and D. S. Scott, {\sl A Compendium of Continuous Lattices},
Springer, Berlin-Heidelberg-New York 1980.

\bibitem{K} Lj. D. R. Kocinac, {\sl Closure properties of function spaces}, Applied General Topology 4 (2) (2003)
255–-261.

\bibitem{LAMPAP} P. TH. Lambrinos and B. K. Papadopoulos,  {\sl The (strong) Isbell topology and (weakly) continuous lattices,}
Continuous Lattices and Applications, Lecture Notes in Pure and Appl. Math.
No. 101, Marcel Dekker, New York 1984, 191--211.

\bibitem{MN} R. McCoy and  I. Ntantu, {\sl Topological properties of spaces of continuous functions}, Lecture Notes in Mathematics, 1315,  Springer Verlang.

\bibitem{N} S. A. Naimpally, {\sl Graph topology for function spaces}, Transactions of the American Mathematical Society, vol. 123, pp. 267–-272, 1966.

\bibitem{XY} Xiaoyong Xi, Jinbo Yang, Coincidence of the Isbell and Scott topologies on domain function spaces, Topology Appl. 164 (2014), 197–-206.

\end{thebibliography}
\end{document}